\documentclass{article}

\usepackage{PRIMEarxiv}

\usepackage[utf8]{inputenc} 
\usepackage[T1]{fontenc}    
\usepackage{hyperref}       
\usepackage{url}            
\usepackage{booktabs}       
\usepackage{amsfonts}       
\usepackage{nicefrac}       
\usepackage{microtype}      
\usepackage{lipsum}
\usepackage{fancyhdr}       
\usepackage{graphicx}       
\graphicspath{{media/}}     

\usepackage{amsmath}
\usepackage{xcolor}
\usepackage{amsthm}
\usepackage{pdflscape}
\usepackage{pgfplots}
\usepackage{mathrsfs}

\usepackage{graphicx}
\usepackage{caption}
\usepackage{subcaption}

\usepackage{booktabs} 
\usepackage{multirow}

\newtheorem{thm}{Theorem}[section]

\theoremstyle{definition}
\newtheorem{defn}[thm]{Definition}

\theoremstyle{remark}

\title{Assessing and Enhancing Graph Neural Networks for Combinatorial Optimization: Novel Approaches and Application in Maximum Independent Set Problems}

\author{
  Chenchuhui Hu \\
  Brandeis University \\
  Waltham\\
  \texttt{chenchuhuih@brandeis.edu} \\
}

\begin{document}
\maketitle

\begin{abstract}
Combinatorial optimization (CO) problems are challenging as the computation time grows exponentially with the input. Graph Neural Networks (GNNs) show promise for researchers in solving CO problems. This study investigates the effectiveness of GNNs in solving the maximum independent set (MIS) problem, inspired by the intriguing findings of Schuetz et al.\cite{schuetz2022combinatorial}, and aimed to enhance this solver. Despite the promise shown by GNNs, some researchers observed discrepancies when reproducing the findings, particularly compared to the greedy algorithm, for instance. We reproduced Schuetz' Quadratic Unconstrained Binary Optimization (QUBO) unsupervised approach and explored the possibility of combining it with a supervised learning approach for solving MIS problems. While the QUBO unsupervised approach did not guarantee maximal or optimal solutions, it provided a solid first guess for post-processing techniques like greedy decoding or tree-based methods. Moreover, our findings indicated that the supervised approach could further refine the QUBO unsupervised solver, as the learned model assigned meaningful probabilities for each node as initial node features, which could then be improved with the QUBO unsupervised approach. Thus, GNNs offer a valuable method for solving CO problems by integrating learned graph structures rather than relying solely on traditional heuristic functions. This research highlights the potential of GNNs to boost solver performance by leveraging ground truth during training and using optimization functions to learn structural graph information, marking a pioneering step towards improving prediction accuracy in a non-autoregressive manner.
\end{abstract}


\section{Introduction}
CO problems seek optimal or suboptimal solutions within a discrete space and are known for their computational complexity. Some well-known CO problems include the Traveling Salesman Problem (TSP), Graph Partition, Graph Similarity, Minimum Vertex Cover (MVC), Graph Coloring (GC), and Maximum Independent Set (MIS). These problems have practical applications in areas such as vehicle routing, scheduling, and calculating Shannon Capacity, particularly using MIS. The application of Shannon Capacity using MIS will be demonstrated in the Application Section.

Traditionally, solvers for CO problems are categorized into three types: exact algorithms, approximation algorithms, and heuristic algorithms. Exact algorithms generally guarantee optimal solutions but are computationally intensive for large-scale problems. Approximation algorithms aim for near-optimal solutions with reduced computational time but may not always be feasible for certain graph-based CO problems such as TSP. Heuristic algorithms usually find suboptimal solutions quickly and efficiently but are less flexible due to each heuristic function being tailored to individual problems \cite{peng2021graph}.

Recent developments in GNNs have shown state-of-the-art performance on various tasks such as node classification, link prediction, and graph classification. GNNs analyze graph data by aggregating structural information, providing a novel approach to solving CO problems. Unlike traditional heuristic functions, which require extensive manual tuning, GNNs offer a more automated solution. This potential has drawn researchers' attention to utilizing GNNs for CO problems, as traditional approaches have limitations. Existing research has explored the potential of GNNs and demonstrated their effectiveness in solving CO problems \cite{schuetz2022combinatorial, physics, sun2023difusco, joshi2019efficient, li2018combinatorial}. However, some studies have reported that GNNs are not yet competitive with standard algorithms like the greedy algorithm \cite{worse_than_greedy, boettcher2023inability}. This discrepancy raised questions about the contributions of GNNs to solving CO problems. It is important to recognize that GNN research in CO problems is still in its early stages compared to the decades of work on traditional methods. We hypothesized that GNNs, by learning the structural information of graphs, could enhance traditional approaches through their inherent self-learning capabilities.

Most CO problems involve graph analysis, generally following two steps: graph representation learning and solving the CO problems using these representations, either with traditional methods or other machine learning approaches \cite{peng2021graph}. The first step can be independent of the second step using typical embedding methods \cite{grover2016node2vec, perozzi2014deepwalk, tang2015line}. However, modern GNNs combine both steps into a unified process where the loss from the downstream task supervises the training of the embeddings \cite{kipf2016semi, hamilton2017inductive, velivckovic2017graph}. Additionally, solvers in the second step can be categorized into non-autoregressive and autoregressive methods. Non-autoregressive methods predict solutions in a single step, allowing for faster inference, while autoregressive methods iteratively build upon partial solutions \cite{peng2021graph}.

Inspired by the findings of Schuetz et al. \cite{schuetz2022combinatorial}, the present research aimed to further investigate the effectiveness of GNNs in solving CO problems, specifically the MIS problem, and to improve the proposed solver to avoid invalid or poor solutions as much as possible. It has been demonstrated that the GNNs presented in Schuetz's paper enhance the final greedy decoding, achieving better results than the simple greedy algorithm. Additionally, we improved the methods from the original paper with better node feature initialization and an optimized QUBO function. Recognizing that nodes with fewer edges are more likely to be part of the independent set, we initialized node features based on their degrees, unlike other methods that use random initialization or set all values to one. The QUBO function was also improved to ensure that the reward for nodes selected in the independent set varied based on their degrees, giving greater rewards to nodes with fewer edges. These improvements allow for strong predictions for the final decoding strategy, effectively solving MIS problems for smaller graphs.

Inspired by numerous studies that have successfully implemented supervised learning approaches \cite{sun2023difusco, li2018combinatorial}, we believed that a fully trained model could effectively learn graph structures and generate meaningful node embeddings for unseen but similar graphs. To verify this, the model was trained using a supervised approach with ground truth generated by the GUROBI optimizer \cite{gurobi2018}, employing a Graph Convolutional Network (GCN) with a normalized adjacency matrix as described by Kipf and Welling \cite{kipf2016semi}. The model architecture was modified to stabilize training, ensuring consistently promising validation results. The QUBO function played a crucial role as an optimization function, aiding in the final prediction for each node. While this remained a non-autoregressive solver, it introduced the concept of leveraging ground truth during training, helping the model learn the structural information of the graph and make initial predictive attempts.

We demonstrated the application of our improved solvers in solving Shannon Channel Capacity, showcasing the application of MIS in the field of information theory. The MIS solver can provide an ideal approximation of Shannon capacity, further detailed in the application section Application \ref{sec:application}.

In summary, this research demonstrated that GNNs can enhance traditional methods, such as the greedy algorithm or degree-based greedy algorithm, in solving CO problems. With the improved node feature initialization and optimized QUBO function, the study shows that GNN approaches like the unsupervised QUBO method can provide stronger initial guesses and more accurate final predictions. Additionally, the trained supervised model can provide better initial node features than degree-based features, making it more effective for MIS problems. Furthermore, we applied these enhanced solvers to calculate Shannon Channel Capacity, highlighting the potential of GNNs to improve prediction accuracy and efficiency in practical applications.

\section{Preliminaries}
\label{sec:prelim}
\subsection{Combinatorial Optimization}
Combinatorial optimization (CO) is a class of the optimization problems in mathematics characterized by discrete solutions. It optimizes the objective function by selecting the optimal subset from a finite set. Finding the optimal solution in CO problem is especially challenging in computer science because it is a NP-complete problem, meaning the complexity grows exponentially with the size of the input space. 

\subsection{Graph Theory}

\begin{defn}
    A \textit{graph} $G$ is a pair $(V(G),E(G))$ where $V(G)$ is any set, whose elements are called \textit{vertices} or \textit{nodes} and $E(G)\subseteq[V(G)]^2:=\{\{x,y\}:x,y\in V(G), x\neq y\}$, whose elements are called \textit{edges}. When $G$ is clear from context, we drop it from the vertex set: $V(G)=V$, and similarly for $E$.
    
    Each node $v_i\in V$ has a $d$-dimensional feature vector $x_i\in\mathbb R^d$, where each element of the feature vector represents a specific property of the node. If there are no specific properties to assign to a node, the feature vector is typically initialized randomly or set to a vector of ones.
\end{defn}

\begin{defn}~ 
A set $I\subseteq V(G)$ is \emph{independent} if any for $x, y \in I$, $(x, y) \notin E(G)$. The \emph{maximum independent set} (or MIS for short) problem is finding $\max\{|I|: I\subseteq V(G)\text{ independent}\}$. This number is often denoted by the symbol $\alpha(G)$.
\end{defn}

If $I$ is maximum independent set, then it is maximal independent set, in the sense that any vertex from $V(G)\setminus I$ has an incident edge with the other endpoint in $I$. However, the converse is generally not true.

\section{Machine Learning and Graph Neural Networks (GNNs)}

\subsection{Graph Neural Networks}

GNNs are widely used for graph representation learning, where each layer aggregates information from neighboring nodes to update node embeddings, effectively capturing local graph structures. These embeddings are used to solve tasks like node classification, link prediction, and graph classification. In CO problems like Maximum Cut and MIS, which involve binary node classification, GNNs classify nodes as 0 or 1, transforming the problem into a node classification task. The final embeddings or predicted probabilities are used to complete the task, with methods like machine learning-based, simple decoding, or greedy decoding.

\subsection{Quadratic Unconstrained Binary Optimization (QUBO)} 
The Quadratic Unconstrained Binary Optimization(QUBO) is a mathematical formulation used in many optimization problems. The goal is to find the final node classification, which is generally binary (0 or 1), minimize the quadratic objective function to approximate the optimization solutions. 

The general form of a QUBO is expressed as $$\min(x^TQx),$$ where $x=(x_{1}, x_{2}, ..., x_{n})$ is the vector of binary decision variables for each node and the matrix $Q$ is a square matrix that is assumed to be either symmetric or upper triangular and represents the coefficients of the quadratic objective function. 

There is a connection between the quadratic unconstrained binary optimization and the Ising model for ferromagnetism in statistical mechanics. The QUBO model uses the following Hamiltonian function to cluster similar points together based on a cost function, which is similar to the one given by the Ising model for magnetic attraction between particles.

The formulation of quadratic objective function is determined by the goal of each CO problems, but the objective function must be quadratic. For example, the goal of MIS is to maximize the size of independent set while ensuring that no two connected nodes are included in the same set. Therefore, the objective function is designed to penalize scenarios where connected nodes are both assigned a value of 1 (indicating inclusion in the set) and to reward configurations that allow as many nodes as possible to be classified as 1 (included in the independent set). The MIS then can be formulated as $$H_{MIS}=-\sum_{i\in V(G)}x_{i}+P\sum_{(i,j)\in E(G)}x_{i}x_{j},$$ where $P$ is the hyper-parameter that configure degree of penalty to neighboring nodes. This is the quadratic objective function so we can convert this to $x^TQx$ as the approximation loss function that we strive to minimize. 

\subsection{Supervised Learning Approach}
Supervised learning is a popular method for addressing CO problems, where a model is trained under the guidance of labels that represent optimal or sub-optimal solutions. This approach enhances the model's accuracy in predicting solutions and effectively tackling CO challenges.

Several existing studies already demonstrated the state of art performance of CO problems with supervised approach. \cite{sun2023difusco,li2018combinatorial} For example, the recent study \cite{sun2023difusco} was motivated by the diffusion model that widely used in the generative model. The idea is to inject noise in the forward process and remove the noise in the reverse process so that aiming to maximize the likelihood of optimal solutions. Li\cite{li2018combinatorial} were using both GCN and tree search to perform supervised learning and generate multiple optimal solutions for the same graph. 

Despite its impressive performance on relatively smaller graphs, supervised learning faces significant hurdles in generating the necessary ground truth data for larger, more complex graphs. Given the NP-hard nature of these problems, producing optimal solutions for extensive or large graphs remains computationally prohibitive, suggesting that the supervised approach may not be ideal under such conditions.

\subsection{Unsupervised Learning Approach}
The unsupervised approach to CO problems involves training models without predefined labels, focusing instead on minimizing a differentiable loss function specifically tailored to each problem. This method has shown promising results, particularly when applied through neural networks.

Recent advancements in this area include the work of Karalias \cite{erdos}, who demonstrated that the probability distribution can be effectively learnt by GNN  with a differentiable probablistic penalty loss to obtain high quality solutions. Schuetz \cite{schuetz2022combinatorial, physics} proposed an unsupervised training process that the loss function is to promote binary variables \{0, 1\} to continuous probability parameters so that loss function can be differentiable. The present research aims to build upon and extend Schuetz's methodology by improving the approach, demonstrating the effectiveness of GNNs, and exploring its potential applications and limitations in more complex scenarios.

\subsection{Debates on the Effectiveness of GNNs for CO Problems}
Following the publication of Schuetz's influential paper \cite{schuetz2022combinatorial}, subsequent research paper raised concern regarding the competitiveness of GNNs in addressing CO problems. Critics argued that GNNs do not surpass or underperform compared to traditional approaches such as greedy algorithm (GA) or degree-based greedy algorithm (DGA), not even mention so much more established methods better than DGA. \cite{worse_than_greedy, boettcher2023inability} They believe GNNs approach to tackle CO problems is not promising and requires more validation. One of the primary challenges of this paper is going to investigate the effectiveness of GNN as so many research has been conducted to prove GNNs are able to solve CO problems.

\section{Methodology}\label{sec: methodology}

\subsection{Unsupervised QUBO Approach}\label{subsec: qubo}
\subsubsection{Node Embedding Initialization}

Most initialization methods for node embeddings, such as random values or setting all to ones, generally work well but may result in poor initial guesses and become trapped in local minima. To avoid these issues, we proposed the degree-based initialization, considering that nodes with fewer connections are more likely to be selected in the independent set. Let $D$ be a vector representing the number of connections for each node, where $d_i \in D$ is the degree of the node $i$. Each node's degree is normalized using Min-max normalization. In the special case where $d_{\text{min}} = d_{\text{max}}$, each $\tilde{d}_i$ is set to a constant value. To ensure nodes with lower normalized degrees are more likely to be selected in the independent set, the final transformation as follows: $$x_i = \frac{1}{(\tilde{d}_i + 1)^k},$$ where k is a constant that controls the range of node feature initialization. When $k=1$ or $k=2$, the range is $[0.5, 1]$ or $[0.25, 1]$ respectively. Moreover, adding 1 in the denominator prevents division by 0, ensuring that nodes without any connections are definitely included in the independent set.

\subsubsection{Hamiltonian Function Modification}
We further modified the Hamiltonian function of the MIS problem. Instead of assigning constant reward 1 for nodes selected in the independent set, we allowed each node to have varying degrees of reward based on their connectivity. This approach ensures that nodes with fewer connections are more likely to be selected. Using the initial node features computed from the degree-based initialization, we can simply generate the reward for each node as follows:
$$r_i=R x_i,$$ where $R$ is a predefined constant reward factor. In our study, we defined $R$ as follows:
$$R=\frac{P|E|}{|V|^n},$$ where $E$ and $V$ are edge set and vertex set respectively, and $P$ is a hyper-parameter indicating the penalty for connected nodes in the independent set, and $n$ is to control the range of rewards the node being selected. The smaller $n$ yields larger rewards and can have better initialization if the graph is large or random like ER. The large and random graphs are complex and can easily trapped in the local minima when initialized which results extremely poor optimization results. Therefore, the modified Hamiltonian function for MIS problem is  $$H_{\text{MIS}}=-\sum_{i\in V(G)}r_i x_{i}+P\sum_{(i,j)\in E(G)}x_i x_j.$$

\subsubsection{Graph Neural Networks Architecture}
We used the similar approach to the Schuetz's paper\cite{schuetz2022combinatorial}, employing simple two-layer graph convolution layers. The output is then fed into component-wise sigmoid transformation to provide the probability for each node. The hidden size is a hyperparameter predefined before training.

Because parameters are initialized randomly, it is likely that they may be poorly initialized. If the QUBO approximation loss is above zero, indicating poor initial parameters, the parameters are reinitialized. This reinitialization step prevents the model from getting stuck in sub-optimal solutions and ensures a better exploration of the solution space. Additionally, we set a penalty threshold to avoid poor minimization results. If the penalty doesn't fall below a certain value, the process reruns to achieve better performance and avoid local minima traps.

\subsubsection{Integration with Greedy Decoding Strategy}
One major drawback of original implementation is that it may generate inaccurate solutions, resulting in two connected nodes being in the independent set. Our approach already mitigated these cases as much as possible, but invalid solutions still present for larger graphs. Therefore, we implemented the simple greedy decoding with the predictions from GNNs to ensure the final solution is valid.

$$\hat{Y} = c_1 \cdot \hat{Y}_{\text{GNNs}} + c_2 \cdot X,$$ where $c_1$ and $c_2$ are hyperparameters controlling the influence of GNN predictions and degree-based initialization in the final prediction. $\hat{Y}_{\text{GNNs}}$ is the probability prediction of each node generated by neural networks, and $X$ is the vector storing the initial node features. Recall that initial node feature is based on the degree of each node, ranging from $[0.5^k, 1]$. The higher node feature means fewer connections and a high likelihood of being selected in the independent set.

After obtaining $\hat{Y}$, the final scores combining both probability predictions from GNNs and degree-based initialization, we use the greedy decoding strategy to find the final solution by sorting $\hat{Y}$. Noting that if $c_1 = 0$, this algorithm reduces to a degree-based greedy algorithm.

\subsection{Supervised + Unsupervised QUBO Approach}

Inspired by previous work that used ground truth to learn the embedding of the graph, we aimed to determine that if the trained model can better grasp the graph structure and further enhance the solver. We generated the ground truth using GUROBI \cite{gurobi2018} solver.

For the supervised learning component, we used the hyperbolic tangent (tanh) activation function to stabilize training, ensuring consistent loss and ensure smooth gradient flow. The tanh function shrinks the range to $[-1, 1]$, keeping outputs after activation functions within a reasonable range. 

The loss function employed is Binary Cross-Entropy (BCE) loss, which effectively handles the binary classification task. The trained model is then used to predict on unseen graphs and output the raw probability as the node feature, and a greedy decoding approach is applied to refine these predictions.

Finally, the results from the supervised learning component can be further improved using the unsupervised QUBO approach, leveraging embeddings from supervised model to enhance overall performance. The whole procedure is exactly the same as in Subsection \ref{subsec: qubo}.

\section{Experiment}\label{sec: experiment}

\subsection{Dataset}

We conducted experiments on various datasets to evaluate the effectiveness of our proposed methods. 

First, random graphs were generated from an Erdős–Rényi (ER) model with 10 nodes (ER10) and 100 nodes (ER100), each consisting of 1000 graphs. This variety ensures the robustness of the model across different graph sizes. Additionally, confusion graphs were generated using 5 vowels with 1-letter, 2-letter, and 3-letter combinations to test the model's performance on structured graphs.

We also utilized 1000 graphs sourced from the SATLIB dataset \cite{hoos2000satlib}, a well-known benchmark for graph-based problems. This dataset was split into 800 graphs for training and validation, and 200 graphs for testing.

\subsection{Metrics}

To evaluate the performance of our methods, we used several metrics, the average size of the independent set (Avg Size), the performance drop from the baseline (Drop) and the time taken to compute the solution (Time).

\subsection{Experimental Setup}

The supervised model was trained for 60 epochs using an RTX 2060 GPU. To ensure the consistency, we set $c1=2$ and $c2=3$ through all experiments. For the unsupervised QUBO approach, we set the epochs to 2000 and for the supervised model, we set the epochs to 200. The penalty threshold was set based on each case.

\begin{table}[h]
    \centering
    \resizebox{\textwidth}{!}{
        \begin{tabular}{l|ccc|ccc|ccc}
            \toprule
            \multirow{2}{*}{Method} & \multicolumn{3}{c|}{ER10} & \multicolumn{3}{c|}{ER100} & \multicolumn{3}{c}{SATLIB} \\
            \cline{2-10}
            & AVG Size & Drop(\%) & Time(s) & AVG Size & Drop(\%) & Time(s) & AVG Size & Drop(\%) & Time(s) \\
            \midrule
            GUROBI* & 3.877 & - & 2.31 & 9.206 & - & 471.69 & 403 & - & 125.12 \\
            DGA & 3.714 & 4.20 & 0.1 & 7.300 & 20.70 & 0.84 & 384.33 & 4.63 & 4.21 \\
            \hline
            QUBO+G & \textbf{3.76} & \textbf{3.01} & 20 & 7.324 & 20.44 & 474 & 386.695 & 4.05 & 73.8 \\
            Supervised+G & 3.733 & 3.71 & 0.2 & 7.333 & 20.35 & 1.6 & 388.695 & 3.55 & 5.87 \\
            Supervised+QUBO+G & 3.756 & 3.12 & 4.59 & \textbf{7.338} & \textbf{20.29} & 62.7 & \textbf{390.51} & \textbf{3.10} & 84 \\
            \bottomrule
        \end{tabular}
    }
    \caption{Comparison of Methods across Different Datasets and Metrics}
    \label{tab:comparison}
\end{table}

\begin{table}[h]
    \centering
    \begin{tabular}{l|ccc|ccc|ccc}
        \toprule
        \multirow{2}{*}{Method} & \multicolumn{3}{c|}{Confusion Graph - 1 letter} & \multicolumn{3}{c|}{Confusion Graph - 2 letter} & \multicolumn{3}{c}{Confusion Graph - 3 letter} \\
        \cline{2-10}
        & Size & Drop(\%) & Time(s) & Size & Drop(\%) & Time(s) & Size & Drop(\%) & Time(s) \\
        \midrule
        GUROBI* & 3.0 & - & 0.135 & 13.0 & - & 0.045 & 63.0 & - & 0.054 \\
        DGA & 3.0 & 0.00 & 0.000 & 11.0 & 15.38 & 0.001 & 46.0 & 26.98 & 0.002 \\
        Supervised+QUBO+G & 3.0 & 0.00 & 0.048 & \textbf{13.0} & \textbf{0.00} & 0.052 & \textbf{63.0} & \textbf{0.00} & 0.063 \\
        \bottomrule

    \end{tabular}

    \caption{Comparison of Traditional DGA and Our Novel Methods across Confusion Graphs}
    \label{tab:Application}
\end{table}
\subsection{Results}
The results are reported in the Table \ref{tab:comparison}. In summary, the supervised model with QUBO further refining has the best performance. Our findings demonstrated that Schultz's approach has significant potential in real applications for enhancing performance and verified the potentials of GNNs.

The findings show that the QUBO unsupervised method generally performs well on regular graphs but struggles with highly random graphs like ER graphs. The randomness in ER graphs makes it difficult for the model to find optimal initial parameters in graph convolutional layers, often leading to poor predictions and local minima trapping. For larger ER graphs, the QUBO method did not significantly outperform traditional degree-based greedy algorithms, and the time spent fine-tuning did not yield proportionate improvements. Simply using the trained model to compute the result is more efficient when considering the time. 

For less random graphs like those in SATLIB, the Supervised+QUBO+G approach outperforms all other methods in our study. Surprisingly, we trained the supervised model with only 800 training samples from SATLIB, yet it proved effective across all types of graph data in our study. Although we did not expect the model to perform well on ER graphs, it still managed to outperform the traditional greedy algorithm without any specific training on ER graphs. Across all experiments and datasets, GNNs with the QUBO Hamiltonian function consistently outperformed traditional degree-based greedy algorithms.

\section{Application to Noisy Communication Channels}\label{sec:application}

This section is dedicated to an application of MIS to information theory. We solved the problem for this particular instances and provide the code in \cite{github}.

Imagine you want to send a message over a noisy channel. That is, there is a risk the message sent is not identical to the one received. For example, elementary school students sometimes struggle to distinguish letters like 'b' and 'd'. So, if a student copies a text, they might swap these two letters by mistake, changing the original message. One way to reduce ambiguity is to add redundancy. In our example, this might include writing the word 'bed' at the beginning of the text. Assuming the mistakes are made consistently, the reader will spot that the letters have been swapped and can then make adjustments before reading the rest of the message. Another solution is to not use confusing letters at all. Either way, a natural question arises in the context of the theory of communication: if I want to communicate with a lexicon of $n$ words effectively, how many letters do I need in my alphabet so that no information is lost? This is a simplified version of what is known as the \emph{Shannon Channel Capacity}.

Now we formalize the setup mathematically. Suppose that $V$ is a set of letters (for example, the 128 ASCII characters, or the 26 letters of the English alphabet). We construct a graph $G$ on $V$ called the \emph{confusion graph} by letting letters be adjacent when there is risk of confusion. In our previous example, $bd$ would be an edge. For the purposes of our code, we generate the edges randomly, since we want our program to solve this problem abstractly and in different circumstances, but naturally, real-life confusion graphs could be used, for instance, for hand-written letters.

To produce a confusion graph for a set of words, we use the strong Cartesian product defined as follows. Given a graph $G=(V,E)$, define $G\times G$ as the graph obtained from the vertex set $V\times V$, where two vertices $(u,x)$ and $(v,y)$ are adjacent if and only if $uv\in E$ and $x=y$, or $u=v$ and $xy\in E$ or $\{uv,xy\}\subseteq E$. In words, two $2$-letter words are at risk of confusion if at least one pair of their letters is at risk of confusion. We iterate this process $k$ times to produce a confusion graph for words of length $k$, and we denote the graph by $G^k$. In symbols, $G^1=G$, $G^2=G\times G$ and $G^{n+1}=G^n\times G$, recursively.

\begin{figure}
\centering
\begin{subfigure}{.5\textwidth}
  \centering
  \includegraphics[width=.8\linewidth]{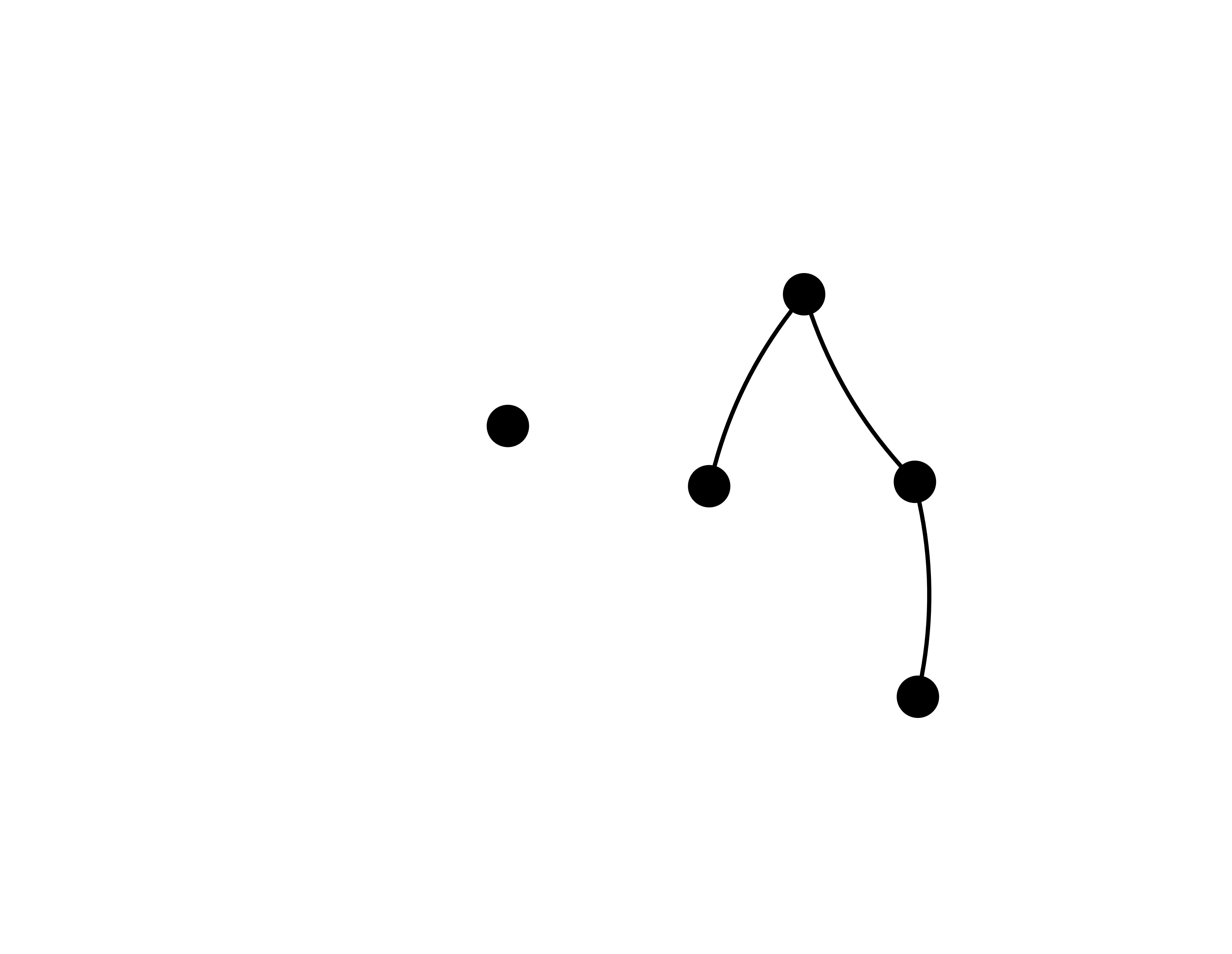}
  \caption{Confusion graph for five vowels}
  \label{fig:sub1}
\end{subfigure}%
\begin{subfigure}{.5\textwidth}
  \centering
  \includegraphics[width=.8\linewidth]{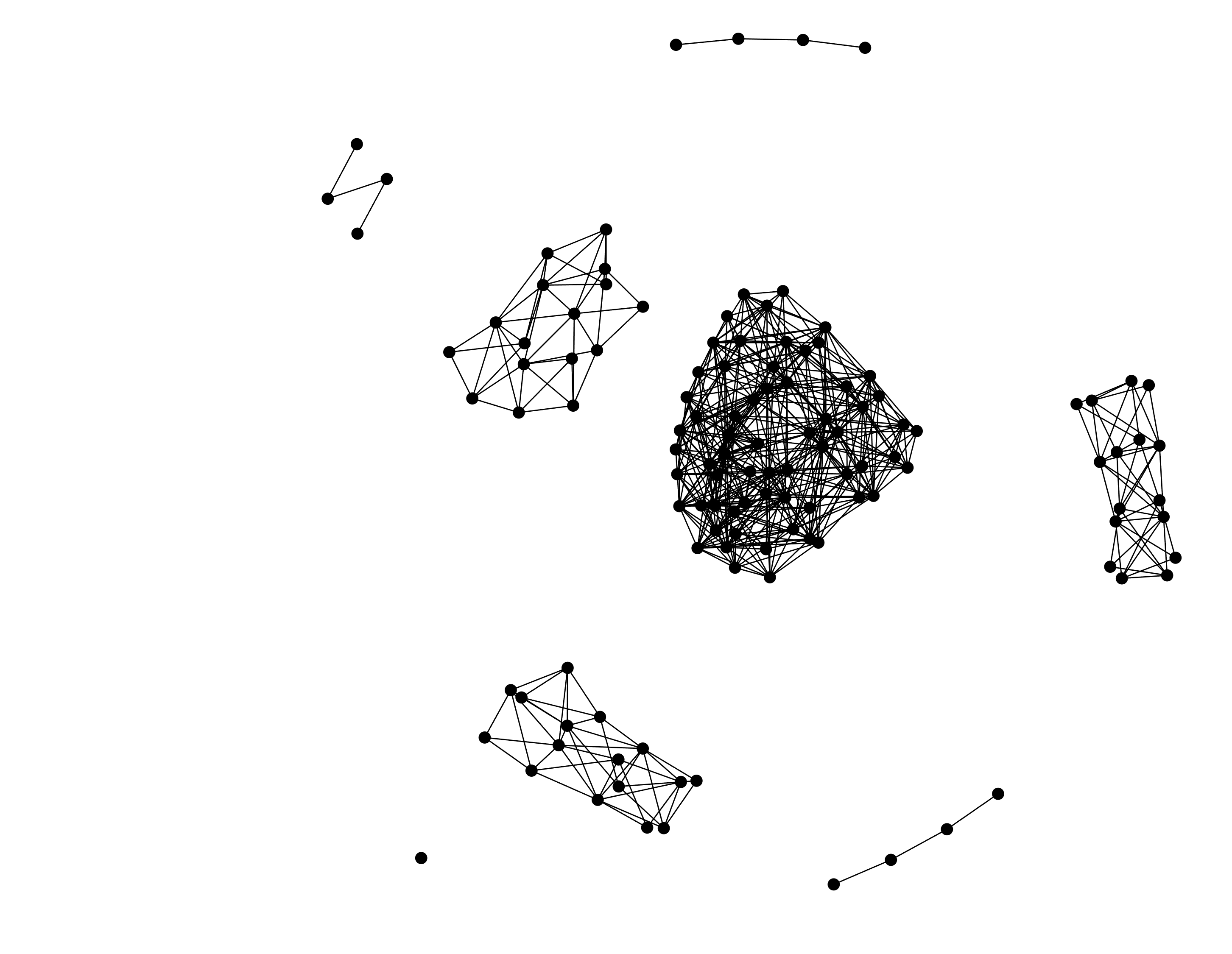}
  \caption{Confusion graph for 3-letter words}
  \label{fig:sub2}
\end{subfigure}
\caption{Confusion graphs by strong product}
\label{fig:test}
\end{figure}

The Shannon capacity (1956) is then $$\lim_{k\to\infty}\sqrt[k]{\alpha(G^k)},$$ where $\alpha(G)$ is a lower bound. For a wide range of graphs, these two numbers are actually equal. See Chapter 14 of \cite{graph}. Thus finding large independent sets in $G^k$ is an interesting problem and one subject to our tools. The size of these graph powers grows rapidly. We illustrate this with an example in Figure~\ref{fig:test}. If our alphabet consists of just the five vowels, the number of 5-letter words already exceeds 3,000.

The GNNs with QUBO solver were implemented to solve the random confusion graph. We made the confusion graphs for five vowels represented in Figure~\ref{fig:test} consisting 1-letter, 2-letter, and 3-letter. The solvers had the ideal approximation of Shannon Channel Capacity. The results were reported in Table\ref{tab:Application}. We compared our results with traditional DGA and showed the superiority of our solver in terms of accuracy and time spent.

\section{Discussion}
The experimental results highlight the effectiveness of GNNs in solving MIS problems across various datasets. After refining the method, Schultz’s approach remains highly effective, particularly in real-world applications.
During the experiments, several notable findings emerged. First, the method performs better on regular graphs than on more random graphs like ER graphs. The QUBO unsupervised approach works well on SATLIB but struggles with ER graphs, likely due to the difficulty of good initialization and the higher chance of getting trapped in local minima. Random graphs often require more fine-tuning and repeated attempts to find optimal or suboptimal solutions, reducing the method’s advantages in time and simplicity.
Second, the supervised model performed surprisingly well on untrained graph types. For instance, although trained on 800 SATLIB graphs, it was effective on ER graphs. We believe the model learned to interpret graph structures and capture essential information, though further validation is needed to explore this generality.
Third, combining the supervised model with the unsupervised QUBO approach produced the best results. This suggests that good node features, often generated by the supervised model, are beneficial for optimization problems, helping the unsupervised model avoid poor local minima.
Despite these advancements, there are limitations. Our study trained on only 800 SATLIB graphs and tested 2,200 in total. For larger ER graphs, this method may not be as efficient in comparison to traditional algorithms. Additionally, our focus was limited to the MIS problem, and further research is needed to apply this approach to other problems. Fine-tuning hyperparameters like penalty thresholds $(c1, c2)$ also proved time-consuming.

\section{Conclusion}
We resolved the problem of invalid solutions and poor approximation, especially less regular graphs, generated from Schultz's QUBO unsupervised approach. We have verified the effectiveness of GNNs on solving combinatorial optimization problems like MIS and can be further enhanced with the supervised model. For this paper, we only experimented with the MIS problems but other experiments may be required to validate if this is applied to other problems. We also provided some insights about the importance of good initialization and avoiding local minima trap, some limitations on the current research, and potential research interest in the future.

\section*{Acknowledgments}
I would like to extend my sincere gratitude to Dr. Tonatiuh Matos Wiederhold for his invaluable support and guidance throughout this project.

\nocite{*}
\bibliographystyle{vancouver}

\newpage
\begin{thebibliography}{10}

\bibitem{schuetz2022combinatorial}
Schuetz MJA, Brubaker JK, Katzgraber HG.
\newblock Combinatorial optimization with physics-inspired graph neural networks.
\newblock Nature Machine Intelligence. 2022;4(4):367-77.

\bibitem{peng2021graph}
Peng Y, Choi B, Xu J.
\newblock Graph learning for combinatorial optimization: a survey of state-of-the-art.
\newblock Data Science and Engineering. 2021;6(2):119-41.

\bibitem{physics}
Schuetz MJA, Brubaker JK, Zhu Z, Katzgraber HG.
\newblock Graph Coloring with Physics-Inspired Graph Neural Networks.
\newblock arXiv preprint arXiv:220201606. 2022.
\newblock Available from: \url{https://arxiv.org/pdf/2202.01606}.

\bibitem{sun2023difusco}
Sun Z, Yang Y.
\newblock {DIFUSCO}: Graph-based Diffusion Solvers for Combinatorial Optimization.
\newblock In: Thirty-seventh Conference on Neural Information Processing Systems; 2023. Available from: \url{https://openreview.net/forum?id=JV8Ff0lgVV}.

\bibitem{joshi2019efficient}
Joshi CK, Laurent T, Bresson X.
\newblock An efficient graph convolutional network technique for the travelling salesman problem.
\newblock arXiv preprint arXiv:190601227. 2019.

\bibitem{li2018combinatorial}
Li Z, Chen Q, Koltun V.
\newblock Combinatorial optimization with graph convolutional networks and guided tree search.
\newblock Advances in neural information processing systems. 2018;31.

\bibitem{worse_than_greedy}
Angelini MC, Ricci-Tersenghi F.
\newblock Modern graph neural networks do worse than classical greedy algorithms in solving combinatorial optimization problems like maximum independent set.
\newblock Nature Machine Intelligence. 2023;5:29-31.

\bibitem{boettcher2023inability}
Boettcher S.
\newblock Inability of a graph neural network heuristic to outperform greedy algorithms in solving combinatorial optimization problems.
\newblock Nature Machine Intelligence. 2023;5(1):24-5.

\bibitem{grover2016node2vec}
Grover A, Leskovec J.
\newblock node2vec: Scalable feature learning for networks.
\newblock In: Proceedings of the 22nd ACM SIGKDD international conference on Knowledge discovery and data mining; 2016. p. 855-64.

\bibitem{perozzi2014deepwalk}
Perozzi B, Al-Rfou R, Skiena S.
\newblock Deepwalk: Online learning of social representations.
\newblock In: Proceedings of the 20th ACM SIGKDD international conference on Knowledge discovery and data mining; 2014. p. 701-10.

\bibitem{tang2015line}
Tang J, Qu M, Wang M, Zhang M, Yan J, Mei Q.
\newblock Line: Large-scale information network embedding.
\newblock In: Proceedings of the 24th international conference on world wide web; 2015. p. 1067-77.

\bibitem{kipf2016semi}
Kipf TN, Welling M.
\newblock Semi-supervised classification with graph convolutional networks.
\newblock arXiv preprint arXiv:160902907. 2016.

\bibitem{hamilton2017inductive}
Hamilton W, Ying Z, Leskovec J.
\newblock Inductive representation learning on large graphs.
\newblock Advances in neural information processing systems. 2017;30.

\bibitem{velivckovic2017graph}
Veli{\v{c}}kovi{\'c} P, Cucurull G, Casanova A, Romero A, Lio P, Bengio Y.
\newblock Graph attention networks.
\newblock arXiv preprint arXiv:171010903. 2017.

\bibitem{gurobi2018}
{Gurobi Optimization, LLC}. Gurobi Optimizer Reference Manual; 2018.
\newblock Available from: \url{http://www.gurobi.com/documentation/8.1/refman/index.html}.

\bibitem{erdos}
Karalias N, Loukas A.
\newblock Erdos Goes Neural: an Unsupervised Learning Framework for Combinatorial Optimization on Graphs.
\newblock Proceedings of NeurIPS. 2020.

\bibitem{hoos2000satlib}
Hoos HH, St{\"u}tzle T.
\newblock SATLIB: An online resource for research on SAT.
\newblock Sat. 2000:283-92.

\bibitem{github}
Matos~Wiederhold T, Hu C. Neuro-Shannon: A Library for Information-Theoretic Learning; 2024.
\newblock Available from: \url{https://github.com/tonamatos/neuro-shannon}.

\bibitem{graph}
Bondy A, Murty USR.
\newblock Graph Theory.
\newblock Springer; 2008.

\bibitem{GNN}
Cappart Q, Ch\'etelat D, Khalil EB, Lodi A, Morris C, Veli{\v{c}}kovi{\'c} P.
\newblock Combinatorial Optimization and Reasoning with Graph Neural Networks.
\newblock arXiv preprint arXiv:210209544. 2022.
\newblock Available from: \url{https://arxiv.org/pdf/2102.09544}.

\bibitem{github1}
Brubaker JK. gcp-with-gnns-example; 2024.
\newblock Version 1.0.0.
\newblock Available from: \url{https://github.com/amazon-science/gcp-with-gnns-example}.

\end{thebibliography}

\end{document}